\pdfoutput=1
\RequirePackage{ifpdf}
\ifpdf % We are running pdfTeX in pdf mode
\documentclass[pdftex]{sigma}
\else
\documentclass{sigma}
\fi

\usepackage{tikz}
\usetikzlibrary{arrows}
\usetikzlibrary{decorations.pathmorphing}
\usetikzlibrary{decorations.markings}
\usetikzlibrary{patterns}
\usetikzlibrary{automata}
\usetikzlibrary{positioning}
\usetikzlibrary{calc}
\usepackage{tikz-cd}
\tikzset{->-/.style={decoration={
			markings,
			mark=at position #1 with {\arrow{latex}}},postaction={decorate}}}

\tikzset{-<-/.style={decoration={
			markings,
			mark=at position #1 with {\arrowreversed{latex}}},postaction={decorate}}}

\usetikzlibrary{shapes.misc}\tikzset{cross/.style={cross out, draw,
		minimum size=2*(#1-\pgflinewidth),
		inner sep=0pt, outer sep=0pt}}

\newcommand{\R}{{\Bbb R}}

\newcommand{\C}{{\Bbb C}}

\DeclareMathOperator{\im}{Im}
\DeclareMathOperator{\re}{Re}

\numberwithin{equation}{section}

\newtheorem{Theorem}{Theorem}[section]

\newtheorem{Lemma}[Theorem]{Lemma}

\newtheorem{assumption}[Theorem]{Assumption}

{ \theoremstyle{definition}
\newtheorem{Definition}[Theorem]{Definition}

\newtheorem{Remark}[Theorem]{Remark} }

\begin{document}
\allowdisplaybreaks

\newcommand{\arXivNumber}{2003.05765}

\renewcommand{\PaperNumber}{079}

\FirstPageHeading

\ShortArticleName{Admissible Boundary Values for the Gerdjikov--Ivanov Equation}

\ArticleName{Admissible Boundary Values for the Gerdjikov--Ivanov\\ Equation with Asymptotically Time-Periodic\\ Boundary Data}

\Author{Samuel FROMM}

\AuthorNameForHeading{S.~Fromm}

\Address{Department of Mathematics, KTH Royal Institute of Technology, 100 44 Stockholm, Sweden}
\Email{\href{mailto:samfro@kth.se}{samfro@kth.se}}

\ArticleDates{Received March 13, 2020, in final form August 09, 2020; Published online August 19, 2020}

\Abstract{We consider the Gerdjikov--Ivanov equation in the quarter plane with Dirichlet boundary data and Neumann value converging to single exponentials $\alpha {\rm e}^{{\rm i}\omega t}$ and $c{\rm e}^{{\rm i}\omega t}$ as $t\to\infty$, respectively.
Under the assumption that the initial data decay as $x\to\infty$, we derive necessary conditions on the parameters $\alpha$, $\omega$, $c$ for the existence of a solution of the corresponding initial boundary value problem.}


\Classification{37K15; 35Q15}

\section{Introduction}

Long time asymptotics of integrable nonlinear partial differential equations (PDEs) can be studied by means of the Riemann--Hilbert (RH) approach. In this approach, which has been successfully applied to several initial value problems on the line, both for decaying and nondecaying initial data, a RH problem is associated to the equation and the asymptotic behavior is computed with the aid of Deift--Zhou nonlinear steepest descent techniques.

\looseness=-1 For initial boundary value problems on the half-line, the RH approach involves additional steps compared to the case on the line, because, in general, not all boundary values are known for a well-posed problem. For instance, if one assumes that the Dirichlet data are given, then the Neumann value has to be computed. This is often referred to as the Dirichlet to Neumann map.

In the case of decaying boundary data, Antonopoulou and Kamvissis~\cite{AntonopoulouKamvissis2015} showed for the defocusing nonlinear Schr\"odinger equation that if the Dirichlet data have sufficient decay as $t\to\infty$, then the Neumann value also decays, thus successfully characterizing the large $t$ limit of the Dirichlet to Neumann map for decaying boundary conditions.

In the setting of nondecaying boundary data, however, less is known. In this paper, we consider the special case of asymptotically periodic boundary values. More specifically, we consider solutions $q(x,t)$ in the quarter plane $\big\{(x,t)\in\R^2\,|\,x\geq 0,\,t\geq 0\big\}$ whose boundary values satisfy
\begin{gather}\label{bddvalues}
q(0,t)\sim \alpha {\rm e}^{{\rm i}\omega t},\qquad q_x(0,t)\sim c {\rm e}^{{\rm i}\omega t},\qquad t\to\infty,
\end{gather}
where $\alpha > 0$, $\omega \in \R$, and $c \in \C$ are three parameters.

For the focusing nonlinear Schr\"odinger equation
\begin{gather}\label{focusingNLS}
{\rm i}q_t+q_{xx}+2|q|^2q=0
\end{gather}
Boutet de Monvel and coauthors \cite{BIK2007,BIK2009,BK2007,BKS2009,MKSZ2010} were able to show that equation \eqref{focusingNLS} has a~solution with boundary values satisfying \eqref{bddvalues} and with decay as $x\to\infty$, if and only if the parameters $(\alpha,\omega,c)$ satisfy either
\begin{gather}\label{admissiblefocNLS1}
c=\pm\alpha\sqrt{\omega-\alpha^2}\qquad\text{and}\qquad\omega\geq\alpha^2
\end{gather}
or
\begin{gather}\label{admissiblefocNLS2}
c=i\alpha\sqrt{|\omega|-2\alpha^2}\qquad\text{and}\qquad\omega\leq-6\alpha^2.
\end{gather}
They also computed the long time asymptotics of any such solution using the Deift--Zhou nonlinear steepest descent method.

The first step in the study of initial boundary value problems whose leading order long-time behaviour is described by a single exponential consists of determining those triples $(\alpha,\omega,c)$ which are admissible.
Here we call a triple $(\alpha,\omega,c)$ admissible if there is a solution of the corresponding initial boundary value problem with boundary values of the form \eqref{focusingNLS} (see Definition \ref{definadmissibility} for the precise definition).
In the case of the focusing NLS equation the admissible parameter triples are precisely those determined by \eqref{admissiblefocNLS1} and \eqref{admissiblefocNLS2}.

The defocusing nonlinear Schr\"odinger equation
\begin{gather*}%\label{defocusingNLS}
{\rm i}q_t+q_{xx}-2|q|^2q=0
\end{gather*}
with boundary values satisfying \eqref{bddvalues} has been studied by Lenells \cite{Ldefocusing-admissible} and Lenells and Fokas \cite{tperiodicI,tperiodicII}. In \cite{Ldefocusing-admissible} it was shown that every admissible parameter triple belongs to one of five families. Note that the corresponding result for the focusing case only leads to two admissible families (cf.\ \eqref{admissiblefocNLS1} and~\eqref{admissiblefocNLS2}). Thus the defocusing case seems to be richer, although it is still unclear if all of the five families determined in~\cite{Ldefocusing-admissible} are indeed admissible.

In this paper we aim to implement the first step in the program initiated by Boutet de Monvel and coauthors described above for the Gerdjikov--Ivanov (GI) equation \cite{Gerdzikov}
\begin{gather}\label{DNLSIII}
{\rm i}q_t+q_{xx}+{\rm i}q^2\bar{q}_x+\frac12 |q|^4q=0.
\end{gather}
Equation \eqref{DNLSIII} is related to the derivative nonlinear Schr\"odinger (DNLS) equation
\begin{gather}\label{DNLS}
{\rm i}u_t+u_{xx}-{\rm i}\big(|u|^2u\big)_x=0
\end{gather}
via the invertible gauge transformation
\begin{gather}\label{gauge}
u(x,t)=q(x,t)\exp\left({\rm i}\int_{x}^\infty|q(y,t)|^2{\rm d}y\right).
\end{gather}

The initial boundary value problem for \eqref{DNLSIII} in the quarter plane $\big\{(x,t)\in\R^2\,|\,x\geq 0,\,t\geq 0\big\}$ is overdetermined in the sense that the Dirichlet and Neumann boundary values at $x = 0$ {\it cannot} both be independently prescribed for a well-posed problem. Indeed, in~\cite{ERDOGANGURELTZIRAKIS2018} it was shown that the Dirichlet initial boundary value problem for \eqref{DNLSIII} is locally well-posed in $H^s([0,\infty))$ for any $s\in\big(\frac12,\frac52\big)$, $s\neq \frac32$, with given initial data $q(x,0)=g(x)$ and Dirichlet boundary data $q(0,t)=h(t)$.
In particular, for any $g\in H^s([0,\infty))$ and $h\in H^{\frac{2s+1}4}([0,\infty))$ satisfying $g(0)=h(0)$, there exists a $T=T\big(\|g\|_{H^s([0,\infty))},\|h\|_{H^{\frac{2s+1}4}([0,\infty))}\big)$ such that this problem has a distributional solution
\begin{gather*}
q\in C_t^0H_x^s([0,T]\times\R)\cap C_x^0H_t^{\frac{2s+1}4}(\R\times[0,T]).
\end{gather*}

We will not give a complete classification of the admissible parameter triples for~\eqref{DNLSIII} but instead focus on two particularly interesting families of parameters.
The first family arises as a~generalization of a two-parameter family of stationary solitons.
The second family arises from the plane wave solutions
\begin{gather*}
q^b(x,t)=\alpha {\rm e}^{{\rm i}\omega t+{\rm i}bx}
\end{gather*}
for suitable parameters $\alpha>0$, $\omega\in\R$, and $b\in\R$.
Within each of these families we give necessary conditions for admissibility.

The proof is inspired by the proof of the corresponding results in \cite{BKS2009} and \cite{Ldefocusing-admissible}.

\section{Main result}

Before stating our main result, we give the definition of an admissible triple (see Definition~1.2 in~\cite{BKS2009} or Definitions~2.1--2.3 in~\cite{Ldefocusing-admissible}) and introduce two special families of parameters.
Let $\mathcal{S}([0,\infty))$ denote the Schwartz space
\begin{gather*}
\mathcal{S}([0,\infty))=\big\{u\in C^\infty([0,\infty))\,|\, \sup_{x\geq 0}\big|x^n u^{(m)}(x)\big|<\infty\text{ for all }n,m=0,1,\dots\big\}.
\end{gather*}

\begin{Definition}\label{definitionsolution}
A solution of the GI equation in the quarter plane $\big\{(x,t)\in\R^2\,|\,x\geq0,\,t\geq 0\big\}$ is a smooth function $q\colon[0,\infty)\times [0,\infty)\to \C$ with $q(\cdot,t)\in \mathcal{S}([0,\infty))$ for each $t\geq 0$, which satisfies~\eqref{DNLSIII} for $x>0$ and $t>0$.
\end{Definition}

\begin{Definition}\label{definadmissibility}
A parameter triple $(\alpha,\omega,c)$ with $\alpha>0$, $\omega\in\R$, and $c\in\C$, is admissible for the GI equation if there exists a solution $q(x,t)$ of~\eqref{DNLSIII} in the quarter plane such that
\begin{gather}\label{asymptoticsadmissible}
q(0,t)-\alpha {\rm e}^{{\rm i}\omega t}\to 0\qquad \text{and}\qquad q_x(0,t)-c {\rm e}^{{\rm i}\omega t}\to 0\text{ sufficiently fast as }t\to\infty.
\end{gather}
\end{Definition}

\begin{Remark} We need a certain order of decay in \eqref{asymptoticsadmissible} to show that certain solutions of Volterra equations are well defined and analytic. For example, it is enough to assume that the order of decay is $O\big(t^{-5/2}\big)$.
\end{Remark}

\subsection{The soliton solution}\label{solitonsolutionintro}

Equation \eqref{DNLS} admits a two-parameter family of solitons \cite{KaupNewell1978} (see also for example (1.2) in~\cite{ColinOhta2006})
\begin{gather*}
u_{\omega,d}(x,t)=\varphi_{\omega,d}(x+dt) \exp \left({\rm i}\omega t- {\rm i}\frac{d}{2}(x+dt) - \frac{3{\rm i}}{4} \int_{x+dt}^{\infty} \varphi_{\omega,d}(y)^{2} {\rm d}y\right),
\end{gather*}
where $d\in\R$, $\omega > d^2/4$, and
\begin{gather*}
\varphi_{\omega,d}(x)=\sqrt{
 \frac{4\omega - d^2}{\omega^{1/2} \left( \cosh \big(\sqrt{4\omega- d^2}x\big)-\frac{d}{2\sqrt{\omega}} \right) } }.
\end{gather*}
Letting $d=0$, applying the gauge transform \eqref{gauge}, and multiplying the resulting function by ${\rm e}^{-{\rm i}\pi/4}$, we obtain a one-parameter family of solutions of the GI equation with periodic boundary values.
More precisely, we obtain that for every $\omega>0$, the function
\begin{gather*}
q_{\omega}(x,t)=\phi_{\omega}(x) {\rm e}^{-\frac{{\rm i}\pi }4+\frac{{\rm i}}{4} \int_x^{\infty } \phi_\omega(y)^2 {\rm d}y}{\rm e}^{{\rm i}\omega t}=\phi_\omega(x){\rm e}^{-{\rm i}\arctan(\tanh(\omega^{1/2}x))}{\rm e}^{{\rm i}\omega t},
\end{gather*}
with
\begin{gather*}
\phi_{\omega}(x) =\sqrt{
 \frac{4\omega}{\omega^{1/2} \cosh \big(\sqrt{4\omega}x\big) } },
\end{gather*}
is a solution of \eqref{DNLSIII} (in the sense of Definition~\ref{definitionsolution}) with boundary values
\begin{gather*}
q_\omega(0,t)=\alpha {\rm e}^{{\rm i}\omega t},\qquad (q_\omega)_x(0,t)=c{\rm e}^{{\rm i}\omega t},
\end{gather*}
where
\begin{gather*}
\alpha=2\omega^{1/4}\qquad\text{and}\qquad c=-2\omega^{3/4}{\rm i}.
\end{gather*}
In particular, it follows that the family of parameters
\begin{gather}\label{solitonfamily}
\big\{\big(\alpha=2\omega^{1/4},\,\omega,\,c=-2\omega^{3/4}{\rm i}\big)\,\big|\,\omega>0\big\}=\left\{\left(\alpha,\,\omega=\frac{\alpha^4}{16},\,c=-\frac{\alpha^3}4{\rm i}\right)\bigg|\, \alpha>0\right\}
\end{gather}
is admissible for the GI equation.

We note that the parameters associated with the soliton solution $q_\omega$ satisfy
\begin{gather}\label{X2=0parameter}
\alpha ^6-2 \alpha ^2 \omega +2|c|^2+4 \alpha ^3 \im(c)=0.
\end{gather}

\subsection{The plane wave}\label{planewaveintro}

Equation \eqref{DNLSIII} admits the plane wave solution
\begin{gather}\label{planewave}
q^b(x,t)=\alpha {\rm e}^{{\rm i}\omega t+{\rm i}bx},
\end{gather}
where $\alpha>0$, $\omega\in\R$, and $b\in\R$ satisfy
\begin{gather}\label{bparameter}
\alpha ^4-2 b^2+2 \alpha ^2 b-2 \omega=0.
\end{gather}
The boundary values of \eqref{planewave} are given by
\begin{gather*}
q^b(0,t)=\alpha {\rm e}^{{\rm i}\omega t},\qquad q^b_x(0,t)=c{\rm e}^{{\rm i}\omega t},
\end{gather*}
where
\begin{gather*}
c=\alpha b {\rm i}.
\end{gather*}
Substituting the latter expression into~\eqref{bparameter}, we find that the parameters associated with the plane wave satisfy the conditions
\begin{align}\label{planewaveparameter}
\re(c)=0,\qquad \im(c)^2+ \alpha^2\omega =\frac{\alpha^6}{2}+\alpha^3 \im(c).
\end{align}
Note that the plane wave \eqref{planewave} itself does not decay as $x\to\infty$ and hence is not a solution of~\eqref{DNLSIII} in the sense of Definition~\ref{definitionsolution}.

\subsection{Statement of the result}

The following theorem classifies all potentially admissible parameter triples within the families corresponding to the stationary soliton and the plane wave given in \eqref{X2=0parameter} and \eqref{planewaveparameter}, respectively.

\begin{Theorem}\label{maintheorem}
Let $\alpha>0$, $\omega\in\R$ and $c\in\C$.
\begin{enumerate}\itemsep=0pt
\item[$(a)$]
Any admissible triple $(\alpha,\omega,c)$ which satisfies \eqref{X2=0parameter} belongs to the family
\begin{gather*}
\left\{\left(\alpha,\,\omega,\,c=\pm{\alpha} \sqrt{ \omega -\frac{\alpha ^4}{16}} -\frac{ \alpha^3}4{\rm i}\right)\!\Bigg|\,\alpha>0,\,\omega\geq\frac{\alpha^4}{16}\right\}\cup\left\{\left(\alpha,\,-\frac{\alpha^4}4, \,-\frac{\alpha^3}2{\rm i}\right)\!\bigg|\,\alpha>0\right\}.
\end{gather*}

\item[$(b)$]
Any admissible triple $(\alpha,\omega,c)$ which satisfies \eqref{planewaveparameter} belongs to one of the two families
\begin{gather*}
\left\{\left(\alpha,\,\omega,\,c=\alpha \left(\frac{\alpha ^2}2-\sqrt{\frac{3 \alpha ^4}4- \omega }\right){\rm i}\right)\Bigg|\,\alpha>0,\,\omega\leq-\frac{\alpha^4}{4}\right\},
\\
\left\{\left(\alpha,\,\omega,\,c=\alpha \left(\frac{\alpha ^2}2+\sqrt{\frac{3 \alpha ^4}4- \omega }\right){\rm i}\right)\Bigg|\,\alpha>0,\,\omega\leq-\frac{1}{2} \big(6 \sqrt{6} +15 \big)\alpha^4\right\}.
\end{gather*}
\end{enumerate}
\end{Theorem}

\begin{Remark}The parameter triples determined in Theorem \ref{maintheorem} are only potentially admissible, i.e., the conditions imposed on a parameter triple $(\alpha,\omega,c)$ by one of the families derived in Theorem \ref{maintheorem} are {necessary} but may not be sufficient for the existence of a solution of \eqref{DNLSIII} with boundary values satisfying \eqref{asymptoticsadmissible}.
It is yet to be determined which of the parameter triples are actually admissible.
In the case of the focusing nonlinear Schr\"odinger equation this was done by constructing an appropriate solution with the help of an associated RH problem \cite{BIK2007,BIK2009,BK2007,BKS2009,MKSZ2010}.
\end{Remark}

\section{Eigenfunctions}

Equation \eqref{DNLSIII} is the compatibility condition of the Lax pair
\begin{gather}\label{laxpair}
\begin{cases}
\mu_x+{\rm i}k^2[\sigma_3,\mu]=U\mu,\\
\mu_t+2{\rm i}k^4[\sigma_3,\mu]=V \mu.
\end{cases}
\end{gather}
Here $k\in\C$ denotes the spectral parameter, $\mu(x,t,k)$ is a $(2\times 2)$-matrix valued eigenfunction and
\begin{gather*}
U=-\frac{{\rm i}}2|q|^2\sigma_3+kQ,
\qquad
\sigma_3=\begin{pmatrix}1&0\\0&-1\end{pmatrix},\qquad Q=\begin{pmatrix}0&q\\ \bar{q}&0\end{pmatrix},
\nonumber
\\
V=-{\rm i}k^2|q|^2\sigma_3+2k^3Q-{\rm i}k Q_x\sigma_3+\frac12(q_x\bar{q}-q\bar{q}_x)\sigma_3+\frac{{\rm i}}4|q|^4\sigma_3.
%\label{definitionV}
\end{gather*}
The above Lax pair arises from the Lax pair for the DNLS equation discovered by Kaup and Newell~\cite{KaupNewell1978} by applying the gauge transformation~\eqref{gauge} (for details see for instance the appendix of~\cite{LiuPerrySulem2016}).
Occasionally it is convenient to consider the rescaled Lax-pair
\begin{align}\label{laxpair-rescaled}
\begin{cases}
\phi_x+{\rm i}k^2\sigma_3\phi=U\phi,
\\
\phi_t+2{\rm i}k^4\sigma_3\phi=V \phi,
\end{cases}
\end{align}
which arises from \eqref{laxpair} through the transformation $\phi=\mu {\rm e}^{-{\rm i}(k^2x+2k^4t)\sigma_3}$.

For the remainder of the paper let $(\alpha,\omega,c)$ be an admissible triple and let $q(x,t)$ be an associated solution of the GI equation in the quarter plane satisfying~\eqref{asymptoticsadmissible}.

\subsection{The background eigenfunction}

Consider the background $t$-part equation
\begin{gather}\label{backgroundtpart}
\phi^b_t+2{\rm i}k^4\sigma_3\phi^b=V^b \phi^b,
\end{gather}
where the matrix $V^b$ is given by $V$ with $q$ and $q_x$ replaced by $\alpha {\rm e}^{{\rm i}\omega t}$ and $c {\rm e}^{{\rm i}\omega t}$, respectively. We define a solution $\phi^b(t,k)$ of \eqref{backgroundtpart} by
\begin{gather*}
\phi^b(t,k)={\rm e}^{\frac{{\rm i}\omega}2t\sigma_3} E(k){\rm e}^{-{\rm i}\Omega(k)t\sigma_3},
\end{gather*}
where $\Omega(k)$ and $E(k)$ are defined by
\begin{gather}\label{Ekdef}
\Omega(k)=\sqrt{4 k^8+2\omega k^4-\frac{\alpha ^6-2 \alpha ^2 \omega +2|c|^2+4 \alpha ^3 \im(c)}2k^2+\frac{\big(\alpha ^4+4 \alpha \im(c)-2 \omega \big)^2}{16}},
\nonumber
\\
E(k)=\sqrt{\frac{2\Omega-H}{2\Omega}}\begin{pmatrix}1&-\dfrac{H}{k\big(\bar{c} + 2{\rm i} \alpha k^2 \big)}\\-\dfrac{H}{k\big({c} - 2{\rm i} \alpha k^2 \big)}&1\end{pmatrix},
\end{gather}
with
\begin{gather*}
H(k)=\Omega(k)-2 k^4+\alpha \im(c)-\alpha ^2 k^2+\frac{\alpha ^4}{4}-\frac{\omega }{2}.
\end{gather*}
We view the functions $\Omega$ and $\sqrt{{(2\Omega-H)}/(2\Omega)}$ as being defined on the cut complex plane $\C\setminus\mathcal{X}_1$ and $\C\setminus\mathcal{X}_2$, respectively, were $\mathcal{X}_i$ contains the branch cuts connecting the zeroes and poles of the respective function.

We have that $\det E(k)=1$ for $k\in\C\setminus\mathcal{X}$ and that $E(k)$ approaches the identity matrix as $k\to\infty$.
Furthermore, the identity
\begin{gather*}
(2\Omega-H)H=-k^2 \big(2 \alpha k^2-{\rm i}\bar{c}\big) \big(2 \alpha k^2+{\rm i}{c}\big)
\end{gather*}
implies that zero is not a branch point of $\sqrt{{(2\Omega-H)}/(2\Omega)}$ so that $E(k)$ is analytic near zero, assuming $0\not\in \mathcal{X}$.

\begin{assumption}\label{branchcutassumption}
We will assume that $\mathcal{X}_1$ and $\mathcal{X}_2$ are invariant under the involutions $k\mapsto -k$ and $k\mapsto \bar{k}$, that $\C\setminus{\mathcal{X}_i}$, $i=1,2$, is connected and that the branch cuts only intersect transversely in at most finitely many points.
\end{assumption}
We will see that in our case the above assumptions are always satisfied.

We fix the branches of $\Omega$ and $\sqrt{{(2\Omega-H)}/(2\Omega)}$ by their asymptotics as $k\to\infty$ as follows:
\begin{gather*}
\Omega(k)=2 k^4+\frac{\omega }{2}+O\big(k^{-2}\big),
\qquad
\sqrt{\frac{2\Omega-H}{2\Omega}}=1+O\big(k^{-2}\big),\qquad k\to\infty.
\end{gather*}
The symmetries of the branch cuts together with the asymptotics of $\Omega$ at infinity imply that $\Omega$ satisfies the identities
\begin{gather*}%\label{Omegasymmetry}
\Omega(k)=\Omega(-k),\qquad k\in\C\setminus\mathcal{X}_1,
\end{gather*}
and
\begin{gather*}%\label{Omegaschwartz}
\Omega(k)=\Omega^*(k),\qquad k\in\C\setminus\mathcal{X}_1,
\end{gather*}
where $\Omega^*(k):=\overline{\Omega(\bar{k})}$ denotes the Schwartz conjugate of $\Omega(k)$.
Similar identities are valid for $\sqrt{{(2\Omega-H)}/{(2\Omega)}}$ on $\C\setminus\mathcal{X}_2$.
In particular, we find
\begin{gather*}
\sigma_1 E(k)^*\sigma_1=E(k),\qquad k\in\C\setminus\mathcal{X},
\end{gather*}
where $\mathcal{X}=\mathcal{X}_1\cup\mathcal{X}_2$ and
\begin{gather*}
\sigma_1=\begin{pmatrix}0&1\\1&0\end{pmatrix}.
\end{gather*}

\subsection{Eigenfunctions}

We define an action $\hat\sigma_3$ on a $2\times 2$ matrix $A$ by $\hat\sigma_3 A=[\sigma_3,A]$, so that ${\rm e}^{\hat \sigma_3}A={\rm e}^{\sigma_3}A{\rm e}^{-\sigma_3}$.
We further define three solutions $\{\phi_j(x,t,j)\}_{j=1}^3$ of~\eqref{laxpair-rescaled} by
\begin{gather*}
\phi_1(x,t,k)=\mu_1(x,t,k){\rm e}^{-{\rm i}(k^2x+(\Omega(k)-\frac\omega2)t)\sigma_3},
\\
\phi_j(x,t,k)=\mu_j(x,t,k){\rm e}^{-{\rm i}(k^2x+2k^4t)\sigma_3},\qquad j=2,3,
\end{gather*}
where $\mu_j$ are $(2\times2)$-matrix valued solutions of the Volterra integral equations
\begin{gather}
\mu_1(x,t,k)={\rm e}^{-{\rm i}k^2x\hat\sigma_3}\bigg\{\mathcal{E}(t,k)-\mathcal{E}(t,k)\int_t^\infty {\rm e}^{{\rm i}(\Omega(k)-\frac\omega2)(t'-t)\hat\sigma_3}\big[\mathcal{E}^{-1}(t',k)
\nonumber
\\
\hphantom{\mu_1(x,t,k)=}{}\times \big(V-V^b\big)(0,t',k)\mu_1(0,t',k)\big]{\rm d}t'+\int_0^x {\rm e}^{{\rm i}k^2x'\hat\sigma_3}[U(x',t)\mu_1(x',t,k)]{\rm d}x'
\bigg\},
\nonumber
\\
\mu_j(x,t,k)=I+\int_{(x_j,t_j)}^{(x,t)}{\rm e}^{{\rm i}[k^2(x'-x)+2k^4(t'-t)]\hat\sigma_3}W_j(x',t',k) ,\qquad j=2,3,\label{definitionmuj}
\end{gather}
with $(x_2,t_2)=(0,0)$, $(x_3,t_3)=(\infty,t)$, and
\begin{gather*}
\mathcal{E}(t,k)={\rm e}^{\frac{{\rm i}\omega}2t\hat\sigma_3}E(k),\qquad W_j=(U{\rm d}x+V{\rm d}t)\mu_j,\quad j=2,3.
\end{gather*}
Finally, we define domains $D_j\subseteq \C$, $j=1,2,3,4$, by
\begin{alignat*}{3}
& D_1=\big\{k\in\C\,|\, \im k^2>0,\, \im \Omega(k)>0\big\},\qquad && D_2=\big\{k\in\C\,|\, \im k^2>0,\,\im \Omega(k)<0\big\},&\\
&D_3=\big\{k\in\C\,|\, \im k^2<0,\,\im \Omega(k)>0\big\},\qquad &&D_4=\big\{k\in\C\,|\, \im k^2<0,\,\im \Omega(k)<0\big\},
\end{alignat*}
and let $D_+=D_1\cup D_3$ and $D_-=D_2\cup D_4$.

Next we will collect some properties of the eigenfunctions $\{\mu_j(x,t,k)\}_{j=1}^3$:
\begin{itemize}\itemsep=0pt
\item The first (resp.\ second) column of $\mu_1(0,t,k)$ is defined and analytic for $k\in D_-\setminus\mathcal{X}$ (resp.\ $D_+\setminus\mathcal{X}$).
Furthermore, the second column of $\mu_1$ has a continuous extension to the boundary of $D_+\setminus\mathcal{X}$, in the sense that away from the branch points the limits from the right and left onto every branch cut in $D_+$ and onto each part of the boundary of $D_+$ exist and are continuous.
Note that if a branch cut can be approached from both right and left from within $D_+\setminus\mathcal{X}$, then the right and left limits are, in general, different.

\item $\mu_2(x,t,k)$ is defined and analytic for all $k\in\C$.

\item The first (resp.\ second) column of $\mu_3(x,t,k)$ is defined and analytic for $\im k^2<0$ (resp.\ $\im k^2>0$) with a continuous extension to $\im k^2\leq 0$ (resp.\ $\im k^2\geq 0$).
\item The $\mu_j$'s are normalized so that
\begin{gather*}
\lim_{t\to\infty} [\mu_1(0,t,k)-\mathcal{E}(t,k)])=0,\qquad k\in(D_-\setminus\mathcal{X},D_+\setminus\mathcal{X}),
\\
\mu_2(0,0,k)=I,\qquad k\in\C,
\\
\lim_{x\to\infty}\mu_3(x,0,k)=I,\qquad k\in\big(\big\{\im k^2\leq 0\big\},\big\{\im k^2\geq0\big\}\big),
\end{gather*}
where $k\in (A_1,A_2)$ indicates that the first and second columns are valid for $k\in A_1$ and $k\in A_2$, respectively.
\end{itemize}

\begin{proof}
The proof is standard, see for instance \cite{DeiftTrubowitz79} or \cite[Proposition 2.2]{tperiodicI}.
The key argument in the proof can be summarized as follows. The first (resp. second) column of the integrand under the $t$-integral appearing in \eqref{definitionmuj} contains the exponential
\begin{gather*}
{\rm e}^{-{\rm i}\Omega(k)(t'-t)}\qquad \big(\text{resp. }{\rm e}^{{\rm i}\Omega(k)(t'-t)}\big),
\end{gather*}
which is bounded in $D_-\setminus\mathcal{X}$ (resp.\ $D_+\setminus\mathcal{X}$).
Furthermore, by assumption~\eqref{asymptoticsadmissible} the term $\big(V-V^b\big)(0,t',k)$ decays as $t\to\infty$.
Standard arguments for Volterra integral equations now imply that the first (resp.\ second) column of $\mu_1(0,t,k)$ is defined and analytic for $k\in D_-\setminus\mathcal{X}$ (resp.\ $D_+\setminus\mathcal{X}$).
The remaining statements follow in a similar fashion.
\end{proof}

\section{Spectral functions}

We define the spectral functions $s(k)$ and $S(k)$ by
\begin{gather*}
s(k) = \mu_3(0,0,k)= \phi_3(0,0,k), \qquad S(k) = \mu_1(0,0,k)= \phi_1(0,0,k).
\end{gather*}
In view of the identities $\sigma_1{\mu_j^*}\sigma_1=\mu_j$, $j=1,2,3,$ we may write
\begin{gather*}
s(k) = \begin{pmatrix} \overline{a(\bar{k})} & b(k) \\
 \overline{b(\bar{k})} & a(k) \end{pmatrix}, \qquad
S(k) = \begin{pmatrix} \overline{A(\bar{k})} & B(k) \\
 \overline{B(\bar{k})} & A(k) \end{pmatrix}.
\end{gather*}
Then
\begin{gather*}
\phi_3(x,t,k) = \phi_2(x,t,k) s(k), \qquad
\phi_1(x,t,k) = \phi_2(x,t,k) S(k).
\end{gather*}

Note that the analyticity properties of $\mu_1$ and $\mu_2$ carry over to $s$ and $S$ and thus to~$a$,~$b$,~$A$ and~$B$.
In particular, the functions $A$ and $B$ are defined and analytic in $D_+\setminus\mathcal{X}$ with a continuous extension to $\bar{D}_+\setminus\mathcal{X}$. Furthermore, away from the branch cuts they also have continuous extensions onto any branch cut intersecting $\bar{D}_+$. The functions $a$ and $b$ are defined and analytic in $\im k^2>0$ with a~continuous extension to $\im k^2\geq 0$.

\section{Global relation}

Consider the $(12)$ entry of the equation
\begin{gather*}
S^{-1}(k)s(k) =
\phi_1^{-1}(0,T,k)\phi_3(0,T,k)
\\
\hphantom{S^{-1}(k)s(k)}{} ={\rm e}^{{\rm i}(\Omega(k)-\frac\omega2)T\sigma_3}\big(\mu_1^{-1}(0,T,k)\mu_3(0,T,k)\big){\rm e}^{-2{\rm i}k^4T\sigma_3}, \qquad k\in D_1\setminus\mathcal{X}.
\end{gather*}
Using the decay of ${\rm e}^{{\rm i}(\Omega(k)+2k^4)T}$, we find
\begin{gather*}
A(k)b(k)-a(k)B(k)=0,\qquad k \in D_1\setminus\mathcal{X},\qquad \im\big(\Omega(k)+2k^4\big)>0.
\end{gather*}
In any unbounded connected component of $D_1\setminus\mathcal{X}$, we can remove the condition $\im\big(\Omega(k)+2k^4\big)>0$ by analytic continuation. Letting~$\mathcal{D}_1$ be any unbounded connected component of~$D_1\setminus\mathcal{X}$, this yields the global relation:
\begin{gather}\label{globalrelation}
A(k)b(k)-a(k)B(k)=0,\qquad k\in\bar{\mathcal{D}}_1.
\end{gather}
\section{Inadmissible triples}

The global relation leads to the following lemma, which is the basis for the proof of Theorem~\ref{maintheorem} (see~\cite{BKS2009} and \cite{Ldefocusing-admissible} for the corresponding result for the focusing and defocusing nonlinear Schr\"odinger equation, respectively).

\begin{Lemma}\label{inadmissiblelemma}
Assume that $\mathcal{D}_1$ is an unbounded connected component of $D_1\setminus\mathcal{X}$ and assume that there exists an open set $U\subseteq \bar{\mathcal{D}}_1$ such that one of the four branch cuts connecting the eight zeroes of~$\Omega^2(k)$ intersects~$U$.
Then the triple $(\alpha,\omega,c)$ is inadmissible.
\end{Lemma}

The proof is standard, see for example \cite[Lemma~3.1]{Ldefocusing-admissible} for the proof of the corresponding result in the case of the defocusing nonlinear Schr\"odinger equation. For the convenience of the reader however, we will present it here as well.

\begin{proof}Let $U\subseteq \bar{\mathcal{D}}_1$ be an open set and $C$ be a branch cut of $\Omega^2(k)$ intersecting $U$.
Note that on $C$ we have $\Omega_+=-\Omega_-$, where $\Omega_+$ and $\Omega_-$ denote the limits of $\Omega$ onto $C$ from the left and right, respectively.
Furthermore, since $U\subseteq \bar{\mathcal{D}}_1\subseteq\bar{D}_1$, we also have $\im \Omega_\pm\geq 0$ on $C\cap U$. Hence $\im \Omega_\pm=0$ on $C\cap U$.
Thus we may define functions $(\mu_1(0,t,k))_\pm$ on $C\cap U$ according to~\eqref{definitionmuj} by replacing~$\mathcal{E}(t,k)$ and $\Omega(k)$ with $\mathcal{E}(t,k)_\pm$ and $\Omega(k)_\pm$, respectively.
We further define eigenfunctions $\nu_\pm(t,k)$ by
\begin{gather}\label{defnu}
(\mu_1(0,t,k))_\pm=\nu_{\pm}(t,k)\mathcal{E}_\pm(t,k),\qquad k\in C\cap U.
\end{gather}
In view of \eqref{definitionmuj} it follows that $\nu_\pm$ satisfies the integral equation
\begin{gather}
\nu_\pm(t,k) =I-\int_t^\infty\phi^b(t,k)\big(\phi^b\big)^{-1}(t',k)
 \big(V-V^b\big)(0,t',k)\nu_\pm(t',k) \nonumber\\
\hphantom{\nu_\pm(t,k) =}{} \times\phi^b(t',k)\big(\phi^b\big)^{-1}(t,k){\rm d}t',\label{volterraeqnnu}
\end{gather}
where $\phi^b(t,k)\big(\phi^b\big)^{-1}(t',k)$ and its inverse are entire functions of~$k$.
The latter statement can be verified directly by computation or one may observe that~$V^b$ is polynomial in~$k$.
Assumption~\eqref{asymptoticsadmissible} yields that $V-V^b=O\big(t^{-5/2}\big)$, which implies that the Volterra equation~\eqref{volterraeqnnu} has a unique solution for $k\in C\cap U$.
Thus $\nu_+=\nu_-=:\nu$.

Let us consider the second column of equation \eqref{defnu}, evaluated at $t=0$, which reads as
\begin{gather*}
\begin{pmatrix}
B(k)\\ A(k)
\end{pmatrix}_\pm
=\nu(0,k)\begin{pmatrix}
E_{12}(k)\\ E_{22}(k)
\end{pmatrix}_\pm.
\end{gather*}
If we write $\nu(0,k)=(\nu_{ij}(k))_{i,j=1,2}$ and use the definition \eqref{Ekdef} of $E(k)$, the last equation can be rewritten as
\begin{gather*}
\left(\frac{B(k)}{A(k)}\right)_\pm
=\frac{\nu_{11}H_\pm-k\big(\bar{c} + 2{\rm i} \alpha k^2 \big)\nu_{12}}{\nu_{21}H_\pm-k\big(\bar{c} + 2{\rm i} \alpha k^2 \big)\nu_{22}}.
\end{gather*}
Using that $H_+-H_-=\Omega_+-\Omega_-=2\Omega_+$ on $C$ and $\det \nu=1$ (which follows from $\det E(k)=1$), we find that
\begin{gather*}
\left(\frac{B(k)}{A(k)}\right)_+-\left(\frac{B(k)}{A(k)}\right)_-
=\frac{-2k\big(\bar{c} + 2{\rm i} \alpha k^2 \big)\Omega_+}{\big(\nu_{21}H_--k\big(\bar{c} + 2{\rm i} \alpha k^2 \big)\nu_{22}\big)\big(\nu_{21}H_+-k\big(\bar{c} + 2{\rm i} \alpha k^2 \big)\nu_{22}\big)}\neq 0.
\end{gather*}
Thus the quotient $B(k)/A(k)$ is discontinuous across $C\cap U$.
Since $a(k)$ and $b(k)$ are continuous in $\bar{D}_1\cup \bar{D}_2$, this contradicts the global relation \eqref{globalrelation}. Hence the triple $(\alpha,\omega,c)$ is inadmissible.
\end{proof}

\section{Proof of Theorem \ref{maintheorem}}

Lemma \ref{inadmissiblelemma} enables us to perform a classification of potentially admissible parameter families.
We do not perform a complete classification as has been done in \cite{BKS2009} and \cite{Ldefocusing-admissible} for the focusing and defocusing nonlinear Schr\"odinger equation, respectively, but instead focus our attention on the two parameter ranges introduced in Sections~\ref{solitonsolutionintro} and~\ref{planewaveintro}.

We note that in the cases below one can directly verify that Assumption~\ref{branchcutassumption} is satisfied by choosing the branch cuts appropriately.

\subsection{The soliton solution case}\label{solitonsection}

In the following we assume that the triple $(\alpha,\omega,c)$ satisfies \eqref{X2=0parameter}.
We write
\begin{gather*}
\Omega(k)=\sqrt{4k^8+X_1 k^4+X_2 k^2+X_3},
\end{gather*}
where
\begin{gather*}
X_1=2\omega,
\qquad X_2=-\frac{\alpha ^6-2 \alpha ^2 \omega +2|c|^2+4 \alpha ^3 \im(c)}{2},
\qquad X_3=\frac{\big(\alpha ^4+4 \alpha \im(c)-2 \omega \big)^2}{16}.
\end{gather*}
Then condition \eqref{X2=0parameter} is equivalent to $X_2=0$.
This implies that
\begin{gather*}
\Omega^2(k) =4\big(k^4-\kappa_+\big)\big(k^4-\kappa_-\big),
\end{gather*}
where
\begin{gather*}
\kappa_\pm=\frac{-X_1\pm\sqrt{X_1^2-16X_3}}8.
\end{gather*}
Solving $X_2=0$ for $\omega$ yields
\begin{gather*}
\omega=\frac{\alpha ^6+2|c|^2+4 \alpha ^3 \im(c)}{2 \alpha ^2},
\end{gather*}
so that
\begin{gather*}
X_3=\frac{|c|^4}{4 \alpha ^4}
\end{gather*}
and
\begin{gather*}
X_1^2-16 X_3 =-\alpha \big(\alpha ^3+4 \im(c)\big) \big(\alpha ^4+4 \alpha \im(c)-4 \omega \big)
\\
\hphantom{X_1^2-16 X_3}{} =\frac{\big(\alpha ^3+4 \im(c)\big) \big(4 \re(c)^2+\big(\alpha ^3+2 \im(c)\big)^2\big)}{\alpha }.
\end{gather*}
We make a case analysis according to the signs of $X_1^2-16X_3$ and $X_1$.

\subsubsection[$X_1^2-16X_3=0$]{$\boldsymbol{X_1^2-16X_3=0}$}

In this case
\begin{gather*}
\Omega^2(k)=\frac{1}{16} \big(8 k^4+X_1\big)^2.
\end{gather*}
Thus $\Omega$ has no branch cuts.
This leads to the following families of potentially admissible triples:
\begin{align}\label{X_2^2-16X_3=0,X2=0}
\left\{\left(\alpha,\,\omega,\,c=\pm{\alpha} \sqrt{ \omega -\frac{\alpha ^4}{16}} -\frac{ \alpha^3}4{\rm i}\right)\Bigg|\,\alpha>0,\,\omega\geq\frac{\alpha^4}{16}\right\}\cup\left\{\left(\alpha,\, -\frac{\alpha^4}4,\,-\frac{\alpha^3}2{\rm i}\right)\bigg|\,\alpha>0\right\}.\!\!\!
\end{align}
Note that the family~\eqref{solitonfamily} is a subset of \eqref{X_2^2-16X_3=0,X2=0}, given by the special case $\omega=\frac{\alpha^4}{16}$.

\subsubsection[$X_1^2-16X_3<0$, $X_1> 0$]{$\boldsymbol{X_1^2-16X_3<0}$, $\boldsymbol{X_1> 0}$}

In this case
\begin{gather*}
\kappa_\pm=\frac{-X_1\pm {\rm i}\sqrt{16 X_3-X_1^2}}8,\qquad \re\kappa_\pm=-\frac{X_1}8<0.
\end{gather*}
Thus each of the sectors created by the rays $\big\{r{\rm e}^{2{\rm i}n\pi/8}\,|\,r\geq 0\big\}$, $n=0,1,2,3,4,5,6,7$, contains exactly one of the eight zeroes~$\Omega^2$.
By using that
\begin{align*}
\{\im\Omega(k)=0\}=\big\{\im\Omega^2(k)=0\big\}\cap\big\{\re\Omega^2(k)\geq0\big\}
\end{align*}
and by directly computing $\im\Omega^2(k)$ and $\re\Omega^2(k)$, we find that
the contour $\im\Omega(k)=0$, shown in Fig.~\ref{figure1}, is given by the eight rays $\big\{r{\rm e}^{2{\rm i}n\pi/8}\,|\,r\geq 0\big\}$, $n=0,1,2,3,4,5,6,7$, together with four simple curves intersecting the rays $\big\{r{\rm e}^{2{\rm i}n\pi/8}\,|\,r\geq 0\big\}$, $n=1,3,5,7$, in one point and connecting the zeroes in the adjoining sectors.

By choosing the branch cut, the component $\mathcal{D}_1$, and the set $U$ appearing in Lemma~\ref{inadmissiblelemma} as shown in Fig.~\ref{figure1-1}, it follows that all parameter triples in this case are inadmissible by Lemma~\ref{inadmissiblelemma}.
Note that while Fig.~\ref{figure1-1} only shows the branch cut in the first quadrant, the remaining branch cuts are chosen in such a way as to satisfy Assumption~\ref{branchcutassumption}.

\subsubsection[$X_1^2-16X_3>0$, $X_1> 0$]{$\boldsymbol{X_1^2-16X_3>0}$, $\boldsymbol{X_1> 0}$}

\begin{figure}[t]\centering
 \begin{tikzpicture}[xscale=0.26,yscale=0.26]
\draw (-10,-10) -- (10,10);
\draw (-10,10) -- (10,-10);
\draw (0,10) -- (0,-10);
\draw (-10,0) -- (10,0);

\draw (3,5) .. controls (3.2,3.2) .. (5,3);
\draw (3,-5) .. controls (3.2,-3.2) .. (5,-3);
\draw (-3,5) .. controls (-3.2,3.2) .. (-5,3);
\draw (-3,-5) .. controls (-3.2,-3.2) .. (-5,-3);

\fill (3,5) circle (0.15);
\fill (3,-5) circle (0.15);
\fill (-3,5) circle (0.15);
\fill (-3,-5) circle (0.15);

\fill (5,3) circle (0.15);
\fill (5,-3) circle (0.15);
\fill (-5,3) circle (0.15);
\fill (-5,-3) circle (0.15);
\end{tikzpicture}
\qquad\qquad
 \begin{tikzpicture}[xscale=0.26,yscale=0.26]

\draw (-10,-10) -- (-4.5,-4.5);
\draw (-3,-3) -- (3,3);
\draw (4.5,4.5) -- (10,10);

\draw (-10,10) -- (-4.5,4.5);
\draw (-3,3) -- (3,-3);
\draw (4.5,-4.5) -- (10,-10);

\fill (-4.5,4.5) circle (0.15);
\fill (-4.5,-4.5) circle (0.15);
\fill (4.5,-4.5) circle (0.15);
\fill (4.5,4.5) circle (0.15);

\fill (-3,3) circle (0.15);
\fill (-3,-3) circle (0.15);
\fill (3,-3) circle (0.15);
\fill (3,3) circle (0.15);

\draw (0,10) -- (0,-10);
\draw (-10,0) -- (10,0);
\end{tikzpicture}
\caption{The qualitative structure of the contour $\im \Omega(k)=0$ (without branch cuts) in the case $X_1^2-16X_3<0$, $X_1>0$ (left) and $X_1^2-16X_3>0$, $X_1>0$ (right).
 The branch points of $\Omega$ are marked with a dot.}\label{figure1}
\end{figure}
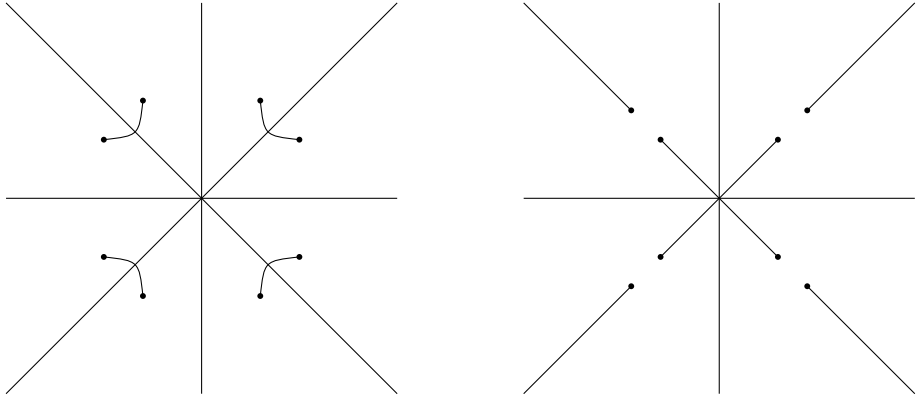

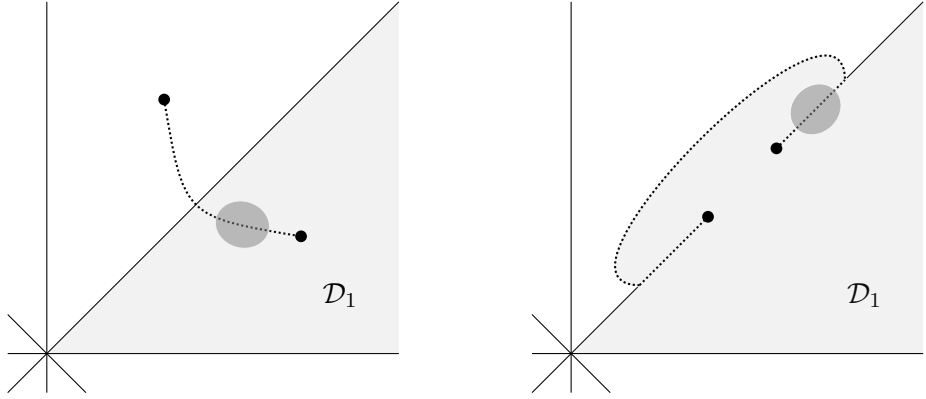
\begin{figure}[t]\centering
 \begin{tikzpicture}[xscale=0.26,yscale=0.26]
\draw (-8,-10) -- (-8,10);
\draw (-10,-8) -- (10,-8);
\draw (-10,-10) -- (10,10);
\draw (-10,-6) -- (-6,-10);

\draw[densely dotted,thick] (-2,5) .. controls (-1,-1) .. (5,-2);

\draw[white,opacity=0.5,fill=gray,rotate around={-15:(2,-1.4)}] (2,-1.4) ellipse (1.4 and 1.2);

\fill[opacity=0.1,fill=gray] (-8,-8)--(10,-8)--(10,10)--cycle;
\node at (7,-5) {$\mathcal{D}_1$};

\fill (-2,5) circle (0.3);
\fill (5,-2) circle (0.3);
\end{tikzpicture}
\qquad\qquad
 \begin{tikzpicture}[xscale=0.26,yscale=0.26]
\draw (-10,-10) -- (-4.6,-4.6);
\draw (6.1,6.1) -- (10,10);

\draw (-10,-6) -- (-6,-10);

\draw[densely dotted,thick] (-4.5,-4.5) to [out=90+90, in=90] (6,6);
\draw[densely dotted,thick] (-4.5,-4.5) -- (-1,-1);
\draw[densely dotted,thick] (2.5,2.5) -- (6,6);
\draw[opacity=0.1,fill=gray] (-4.5,-4.5) to [out=90+90, in=90] (6,6);

\draw[white,opacity=0.5,fill=gray,rotate around={45:(4.5,4.5)}] (4.5,4.5) ellipse (1.4 and 1.2);

\node at (7,-5) {$\mathcal{D}_1$};

\fill[opacity=0.1,fill=gray] (-8,-8)--(10,-8)--(10,10)--cycle;

\fill (2.5,2.5) circle (0.3);
\fill (-1,-1) circle (0.3);

\draw (-8,-10) -- (-8,10);
\draw (-10,-8) -- (10,-8);
\end{tikzpicture}
\caption{A possible choice of branch cuts in the first quadrant in the case $X_1^2-16X_3<0$, $X_1>0$ (left) and $X_1^2-16X_3>0$, $X_1>0$ (right). The branch points of $\Omega$ are marked with a dot and the branch cuts are represented by dotted lines. The set $U$ is shaded dark gray and the set $\mathcal{D}_1$ is shaded light gray.}\label{figure1-1}
\end{figure}

In this case we find that $\kappa_\pm<0$ and that $\kappa_+-\kappa_-=\frac14{\sqrt{X_1^2-16X_3}}>0$.
The contour $\im\Omega(k)=0$, shown in Fig.~\ref{figure1}, is given by the coordinate axes together with the four rays $\big\{r{\rm e}^{2{\rm i}n\pi/4+{\rm i}\pi/4}\,|\,r\geq 0\big\}$, $n=0,1,2,3$, excluding the parts of the rays connecting the four pair of zeroes
\begin{gather*}
|\kappa_\pm|^{1/4}{\rm e}^{2{\rm i}n\pi/4+{\rm i}\pi/4},\qquad n=0,1,2,3.
\end{gather*}
By choosing the branch cut, the component $\mathcal{D}_1$, and the set $U$ as shown in Fig.~\ref{figure1-1}, it follows that all parameter triples in this case are inadmissible by Lemma \ref{inadmissiblelemma}.

\subsubsection[$X_1^2-16X_3>0$, $X_1\leq 0$]{$\boldsymbol{X_1^2-16X_3>0}$, $\boldsymbol{X_1\leq 0}$}

This case is empty.
Indeed, $X_1^2-16 X_3>0$ implies $\im(c)>-\frac{\alpha^3}4$ so that
\begin{align*}
X_1=2\omega=\omega=\frac{\alpha ^6+2|c|^2+4 \alpha ^3 \im(c)}{ \alpha ^2}>\frac{2|c|^2}{\alpha^2}>0.
\end{align*}

\subsubsection[$X_1^2-16X_3<0$, $X_1< 0$]{$\boldsymbol{X_1^2-16X_3<0}$, $\boldsymbol{X_1< 0}$}

In this case
\begin{gather*}
\kappa_\pm=\frac{-X_1\pm {\rm i}\sqrt{16 X_3-X_1^2}}8,\qquad \re\kappa_\pm=-\frac{X_1}8>0.
\end{gather*}
Thus $\arg(\kappa_+)\in(0,\pi/4)$ and $\arg(\kappa_-)\in(-\pi/4,0)$. The corresponding roots of $\Omega^2$ thus have arguments
\begin{gather*}
(0,\pi/16)+\frac{2\pi {\rm i}n}4\qquad \text{and} \qquad(-\pi/16,0)+\frac{2\pi {\rm i}n}4,\qquad n=0,1,2,3.
\end{gather*}
Consequently, each of the sectors created by the eight rays $\big\{r{\rm e}^{2{\rm i}n\pi/8}\,|\,r\geq 0\big\}$, $n=0,1,2,3,4,5,\allowbreak 6,7$, contains exactly one zero of $\Omega^2$.
The contour $\im\Omega(k)=0$, shown in Fig.~\ref{figure2}, is given by the rays $\big\{r{\rm e}^{2{\rm i}n\pi/8}\,|\,r\geq 0\big\}$, $n=0,1,2,3,4,5,6,7$, together with curves intersecting the rays $\big\{r{\rm e}^{2{\rm i}n\pi/8}\,|\,r\geq 0\big\}$, $n=0,2,4,6$, in one point and connecting the zeroes in the adjoining sectors.

By choosing the branch cut, the component $\mathcal{D}_1$, and the set $U$ as shown in Fig.~\ref{figure2-1}, it follows that all parameter triples in this case are inadmissible by Lemma~\ref{inadmissiblelemma}.

\subsubsection[$X_1^2-16X_3<0$, $X_1= 0$]{$\boldsymbol{X_1^2-16X_3<0}$, $\boldsymbol{X_1= 0}$}
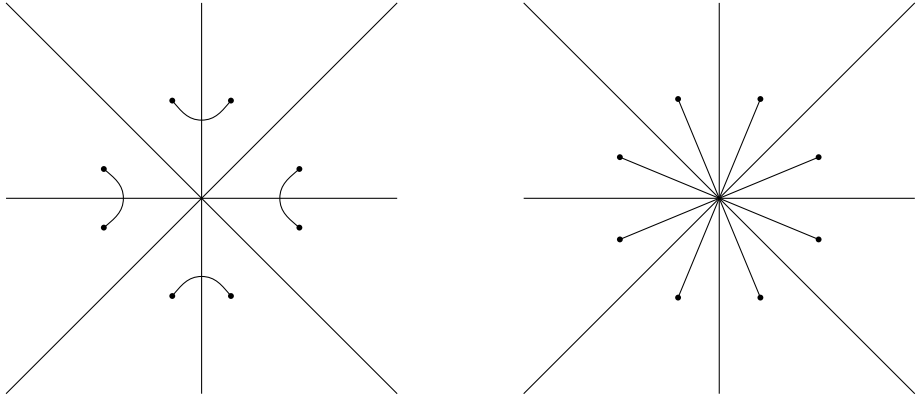
\begin{figure}[t]\centering
 \begin{tikzpicture}[xscale=0.26,yscale=0.26]
\draw (-10,-10) -- (10,10);
\draw (-10,10) -- (10,-10);
\draw (0,10) -- (0,-10);
\draw (-10,0) -- (10,0);

\draw (4,0) to [out=90, in=220] (5,1.5);
\draw (4,0) to [out=270, in=140] (5,-1.5);

\draw (-4,0) to [out=90, in=320] (-5,1.5);
\draw (-4,0) to [out=270, in=40] (-5,-1.5);

\draw (0,4) to [out=0, in=230] (1.5,5);
\draw (0,4) to [out=180, in=310] (-1.5,5);

\draw (0,-4) to [out=0, in=130] (1.5,-5);
\draw (0,-4) to [out=180, in=410] (-1.5,-5);

\fill (5,1.5) circle (0.15);
\fill (-5,1.5) circle (0.15);
\fill (5,-1.5) circle (0.15);
\fill (-5,-1.5) circle (0.15);
\fill (1.5,5) circle (0.15);
\fill (-1.5,5) circle (0.15);
\fill (1.5,-5) circle (0.15);
\fill (-1.5,-5) circle (0.15);
\end{tikzpicture}
\qquad\qquad
 \begin{tikzpicture}[xscale=0.26,yscale=0.26]
\draw (-10,-10) -- (10,10);
\draw (-10,10) -- (10,-10);
\draw (0,10) -- (0,-10);
\draw (-10,0) -- (10,0);

\draw (0,0) -- (5.5*0.92388, 5.5*0.382683);
\draw (0,0) -- (5.5*0.382683,5.5*0.92388);

\fill (5.5*0.92388, 5.5*0.382683) circle (0.15);
\fill (5.5*0.382683,5.5*0.92388) circle (0.15);

\draw (0,0) -- (-5.5*0.92388, 5.5*0.382683);
\draw (0,0) -- (-5.5*0.382683,5.5*0.92388);

\fill (-5.5*0.92388, 5.5*0.382683) circle (0.15);
\fill (-5.5*0.382683,5.5*0.92388) circle (0.15);

\draw (0,0) -- (-5.5*0.92388, -5.5*0.382683);
\draw (0,0) -- (-5.5*0.382683,-5.5*0.92388);

\fill (-5.5*0.92388, -5.5*0.382683) circle (0.15);
\fill (-5.5*0.382683,-5.5*0.92388) circle (0.15);

\draw (0,0) -- (5.5*0.92388, -5.5*0.382683);
\draw (0,0) -- (5.5*0.382683,-5.5*0.92388);

\fill (5.5*0.92388, -5.5*0.382683) circle (0.15);
\fill (5.5*0.382683,-5.5*0.92388) circle (0.15);

\end{tikzpicture}
\caption{The qualitative structure of the contour $\im \Omega(k)=0$ (without branch cuts) in the case $X_1^2-16X_3<0$, $X_1<0$ (left) and $X_1^2-16X_3<0$, $X_1=0$ (right). The branch points of~$\Omega$ are marked with a dot.}\label{figure2}
\end{figure}

This case is equivalent to $X_3>0$, $X_1=X_2=0$.
In this case
\begin{gather*}
\Omega^2(k)=4k^8+X_3,
\end{gather*}
so that the roots of $\Omega^2$ are given by
\begin{gather*}
\left(\frac{X_3}{4}\right)^{1/8}{\rm e}^{2{\rm i}n\pi/8+{\rm i}\pi/8},\qquad n=0,1,2,3,4,5,6,7.
\end{gather*}
The \looseness=1 contour $\im\Omega(k)=0$, shown in Fig.~\ref{figure2}, is given by the eight rays $\big\{r{\rm e}^{2{\rm i}n\pi/8}\,|\,r\geq 0\big\}$, $n=0,1,2,3,4,5,6,7$, together with straight lines connecting the origin with each of the zeroes of~$\Omega^2$.

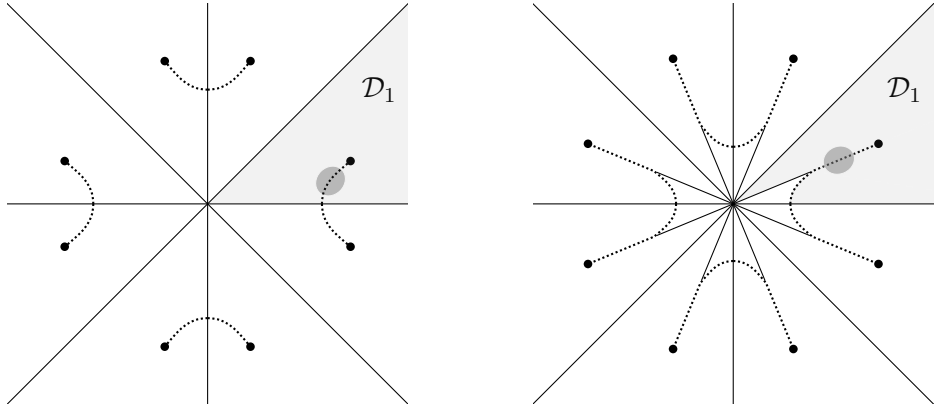
\begin{figure}[t]\centering
 \begin{tikzpicture}[xscale=0.38,yscale=0.38]%[xscale=0.5,yscale=0.5]
 \clip (-7,-7) rectangle (7,7);
\draw (-10,-10) -- (10,10);
\draw (-10,10) -- (10,-10);
\draw (0,10) -- (0,-10);
\draw (-10,0) -- (10,0);

\draw[white,opacity=0.5,fill=gray,rotate around={45:(4.3,0.8)}] (4.3,0.8) ellipse (0.56 and 0.48);

\draw[densely dotted,thick] (4,0) to [out=90, in=220] (5,1.5);
\draw[densely dotted,thick] (4,0) to [out=270, in=140] (5,-1.5);

\draw[densely dotted,thick] (-4,0) to [out=90, in=320] (-5,1.5);
\draw[densely dotted,thick] (-4,0) to [out=270, in=40] (-5,-1.5);

\draw[densely dotted,thick] (0,4) to [out=0, in=230] (1.5,5);
\draw[densely dotted,thick] (0,4) to [out=180, in=310] (-1.5,5);

\draw[densely dotted,thick] (0,-4) to [out=0, in=130] (1.5,-5);
\draw[densely dotted,thick] (0,-4) to [out=180, in=410] (-1.5,-5);

\fill (5,1.5) circle (0.15);
\fill (-5,1.5) circle (0.15);
\fill (5,-1.5) circle (0.15);
\fill (-5,-1.5) circle (0.15);
\fill (1.5,5) circle (0.15);
\fill (-1.5,5) circle (0.15);
\fill (1.5,-5) circle (0.15);
\fill (-1.5,-5) circle (0.15);

\node at (6,4) {$\mathcal{D}_1$};

\fill[opacity=0.1,fill=gray] (-0,-0)--(10,-0)--(10,10)--cycle;

\end{tikzpicture}
\qquad\qquad
 \begin{tikzpicture}[xscale=0.38,yscale=0.38]%[xscale=0.5,yscale=0.5]
 \clip (-7,-7) rectangle (7,7);
\draw (-10,-10) -- (10,10);
\draw (-10,10) -- (10,-10);
\draw (0,10) -- (0,-10);
\draw (-10,0) -- (10,0);

\draw (0,0) -- (2.9*0.92388, 2.9*0.382683);
\draw (0,0) -- (2.9*0.382683,2.9*0.92388);
\draw[densely dotted,thick] (3.05*0.92388, 3.05*0.382683) -- (5.5*0.92388, 5.5*0.382683);
\draw[densely dotted,thick] (3.05*0.382683,3.05*0.92388) -- (5.5*0.382683,5.5*0.92388);

\fill (5.5*0.92388, 5.5*0.382683) circle (0.15);
\fill (5.5*0.382683,5.5*0.92388) circle (0.15);

\draw (0,0) -- (-2.9*0.92388, 2.9*0.382683);
\draw (0,0) -- (-2.9*0.382683,2.9*0.92388);
\draw[densely dotted,thick] (-3.05*0.92388, 3.05*0.382683) -- (-5.5*0.92388, 5.5*0.382683);
\draw[densely dotted,thick] (-3.05*0.382683,3.05*0.92388) -- (-5.5*0.382683,5.5*0.92388);

\fill (-5.5*0.92388, 5.5*0.382683) circle (0.15);
\fill (-5.5*0.382683,5.5*0.92388) circle (0.15);

\draw (0,0) -- (-2.9*0.92388, -2.9*0.382683);
\draw (0,0) -- (-2.9*0.382683,-2.9*0.92388);
\draw[densely dotted,thick] (-3.05*0.92388, -3.05*0.382683) -- (-5.5*0.92388, -5.5*0.382683);
\draw[densely dotted,thick] (-3.05*0.382683,-3.05*0.92388) -- (-5.5*0.382683,-5.5*0.92388);

\fill (-5.5*0.92388, -5.5*0.382683) circle (0.15);
\fill (-5.5*0.382683,-5.5*0.92388) circle (0.15);

\draw (0,0) -- (2.9*0.92388, -2.9*0.382683);
\draw (0,0) -- (2.9*0.382683,-2.9*0.92388);
\draw[densely dotted,thick] (3.05*0.92388, -3.05*0.382683) -- (5.5*0.92388, -5.5*0.382683);
\draw[densely dotted,thick] (3.05*0.382683,-3.05*0.92388) -- (5.5*0.382683,-5.5*0.92388);

\fill (5.5*0.92388, -5.5*0.382683) circle (0.15);
\fill (5.5*0.382683,-5.5*0.92388) circle (0.15);

\node at (6,4) {$\mathcal{D}_1$};

\draw[white,opacity=0.5,fill=gray,rotate around={25:(4*0.92388, 4*0.382683)}] (4*0.92388, 4*0.382683) ellipse (0.56 and 0.48);

\draw[densely dotted,thick] (2,0) to [out=90, in=220] (3*0.92388, 3*0.382683);
\draw[densely dotted,thick] (2,0) to [out=270, in=140] (3*0.92388, -3*0.382683);

\draw[opacity=0.1,fill=gray] (2,0) to [out=90, in=220] (3*0.92388, 3*0.382683);
\fill[opacity=0.1,fill=gray] (3*0.92388, 3*0.382683)--(10, 10.8239*0.382683)--(10,0)--cycle;
\fill[opacity=0.1,fill=gray] (-0,-0)--(10, 10.8239*0.382683)--(10,10)--cycle;
\fill[opacity=0.1,fill=gray] (3*0.92388, 3*0.382683)--(10, 0)--(2,0)--cycle;

\draw[densely dotted,thick] (0,2) to [out=0, in=230] (3*0.382683,3*0.92388);
\draw[densely dotted,thick] (0,2) to [out=180, in=310] (-3*0.382683,3*0.92388);

\draw[densely dotted,thick] (-2,0) to [out=90, in=320] (-3*0.92388, 3*0.382683);
\draw[densely dotted,thick] (-2,0) to [out=270, in=40] (-3*0.92388, -3*0.382683);

\draw[densely dotted,thick] (0,-2) to [out=0, in=130] (3*0.382683,-3*0.92388);
\draw[densely dotted,thick] (0,-2) to [out=180, in=410] (-3*0.382683,-3*0.92388);
\end{tikzpicture}
\caption{A possible choice of branch cuts in the case $X_1^2-16X_3<0$, $X_1<0$ (left) and $X_1^2-16X_3<0$, $X_1=0$ (right).
 The branch points of $\Omega$ are marked with a dot and the branch cuts are represented by dotted lines. The set $U$ is shaded dark gray and the set $\mathcal{D}_1$ is shaded light gray.}\label{figure2-1}
\end{figure}

By choosing the branch cut, the component $\mathcal{D}_1$, and the set~$U$ as shown in Fig.~\ref{figure2-1}, it follows that all parameter triples in this case are inadmissible by Lemma~\ref{inadmissiblelemma}.

\subsection{The plane wave case}\label{planewavesection}

Let $(\alpha,\omega,c)$ belong satisfy~\eqref{planewaveparameter}. We introduce a new parameter $b=-c{\rm i}/\alpha$.
Then we may write
\begin{gather*}
\Omega^2(k)=\frac{1}{4} \big(b-2 k^2\big)^2 \big(\big(b+2 k^2\big)^2+\alpha^2\big(2 b+\alpha^2\big)\big).
\end{gather*}
The zeroes of $\Omega^2(k)$ are given by
\begin{gather*}
\pm\frac{\sqrt{b}}{\sqrt{2}}\text{ (double)},\qquad \frac{-\alpha \pm {\rm i} \sqrt{\alpha^2+2 b}}2, \qquad \frac{\alpha \pm {\rm i} \sqrt{\alpha^2+2 b}}2.
\end{gather*}

\subsubsection[$b\leq -\frac{\alpha^2}2$]{$\boldsymbol{b\leq -\frac{\alpha^2}2}$}
In this case all the branch points of $\Omega$ lie on the real axis.
Thus we cannot rule out the corresponding triples using Lemma~\ref{inadmissiblelemma}.
This leads to the following family of potentially admissible triples
\begin{gather*}
\left\{\left(\alpha,\, \omega=\frac{\alpha ^4}{2}-b^2+\alpha ^2 b,\,c=\alpha b{\rm i}\right)\bigg|\, \alpha^2+2b\leq 0,\, \alpha>0\right\}
\end{gather*}
or
\begin{gather*}
\left\{\left(\alpha,\,\omega,\,c=\alpha \left(\frac{\alpha ^2}2-\sqrt{\frac{3 \alpha ^4}4- \omega }\right){\rm i}\right)\Bigg|\, \alpha>0,\,\omega\leq-\frac{\alpha^4}{4}\right\}.
\end{gather*}

\subsubsection[$-\frac{\alpha^2}2<b<\big(2+\sqrt{6}\big)\alpha^2$]{$\boldsymbol{-\frac{\alpha^2}2<b<\big(2+\sqrt{6}\big)\alpha^2}$}
In this case $\Omega$ has a branch point in each quadrant of the complex plane.
The contour $\im \Omega(k)=0$, shown in Fig.~\ref{figure3}, consists of the coordinate axes together with four simple curves starting from the four branch points $\frac{\alpha \pm {\rm i}\sqrt{\alpha^2+2 b}}2$ and $\frac{-\alpha \pm {\rm i}\sqrt{\alpha^2+2 b}}2$ and asymptoting towards the curves ${\rm e}^{\pm {\rm i}\frac{\pi}4}$ and ${\rm e}^{{\rm i}\pi \mp {\rm i}\frac{\pi}4}$, respectively.

By choosing the branch cut, the component~$\mathcal{D}_1$, and the set~$U$ as shown in Fig.~\ref{figure3-1}, it follows that all parameter triples in this case are inadmissible by Lemma~\ref{inadmissiblelemma}.

\subsubsection[$\big(2+\sqrt{6}\big)\alpha^2\leq b$]{$\boldsymbol{\big(2+\sqrt{6}\big)\alpha^2\leq b}$}
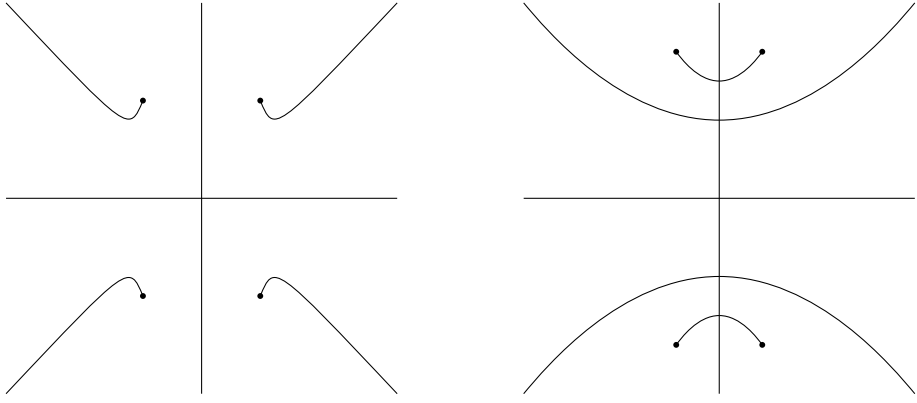
\begin{figure}[t]\centering
 \begin{tikzpicture}[xscale=0.26,yscale=0.26]
\draw (0,10) -- (0,-10);
\draw (-10,0) -- (10,0);

\draw (3,5) .. controls (3.7,3.3) .. (10,10);
\draw (3,-5) .. controls (3.7,-3.3) .. (10,-10);
\draw (-3,5) .. controls (-3.7,3.3) .. (-10,10);
\draw (-3,-5) .. controls (-3.7,-3.3) .. (-10,-10);

\fill (3,5) circle (0.15);
\fill (3,-5) circle (0.15);
\fill (-3,5) circle (0.15);
\fill (-3,-5) circle (0.15);
\end{tikzpicture}
\qquad\qquad
 \begin{tikzpicture}[xscale=0.26,yscale=0.26]
\draw (0,10) -- (0,-10);
\draw (-10,0) -- (10,0);

\draw (0,6) parabola (-2.2,7.5);
\draw (0,6) parabola (2.2,7.5);
\draw (0,4) parabola (-10,10);
\draw (0,4) parabola (10,10);

\fill (-2.2,7.5) circle (0.15);
\fill (2.2,7.5) circle (0.15);
\draw (0,-6) parabola (-2.2,-7.5);
\draw (0,-6) parabola (2.2,-7.5);
\draw (0,-4) parabola (-10,-10);
\draw (0,-4) parabola (10,-10);

\fill (-2.2,-7.5) circle (0.15);
\fill (2.2,-7.5) circle (0.15);
\end{tikzpicture}
\caption{The qualitative structure of the contour $\im \Omega(k)=0$ (without branch cuts) in the case $-\frac{\alpha^2}2<b<\big(2+\sqrt{6}\big)\alpha^2$ (left) and $\big(2+\sqrt{6}\big)\alpha^2\leq b$ (right).
 The branch points of $\Omega$ are marked with a dot.}\label{figure3}
\end{figure}

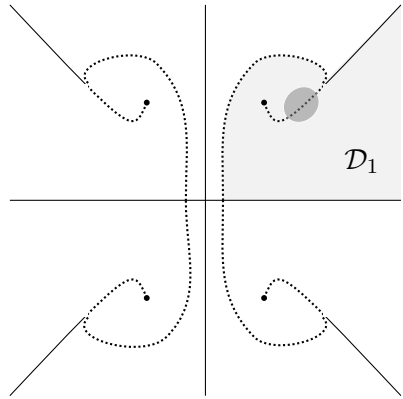
\begin{figure}[t]\centering
 \begin{tikzpicture}[xscale=0.26,yscale=0.26]
\draw (0,10) -- (0,-10);
\draw (-10,0) -- (10,0);

\begin{scope}
\clip (-6,-6) rectangle (6,6);
\draw[densely dotted,thick] (3,5) .. controls (3.7,3.3) .. (10,10);
\draw[densely dotted,thick] (3,-5) .. controls (3.7,-3.3) .. (10,-10);
\draw[densely dotted,thick] (-3,5) .. controls (-3.7,3.3) .. (-10,10);
\draw[densely dotted,thick] (-3,-5) .. controls (-3.7,-3.3) .. (-10,-10);
\end{scope}

 \draw [densely dotted,thick] (6,6) to[out=45,in=45] (2.5,7) to[out=45-180,in=90] (0.9,0);
 \draw [densely dotted,thick] (6,-6) to[out=-45,in=-45] (2.5,-7) to[out=-45+180,in=-90] (0.9,0);

 \draw [densely dotted,thick] (-6,6) to[out=45+90,in=45+90] (-2.5,7) to[out=45-180+90,in=90] (-1,0);
 \draw [densely dotted,thick] (-6,-6) to[out=-45-90,in=-45-90] (-2,-7) to[out=-45+180-90,in=-90] (-1,0);

 \fill[opacity=0.1,fill=gray] (6,6) to[out=45,in=45] (2.5,7) to[out=45-180,in=90] (0.9,0);
 \fill[opacity=0.1,fill=gray] (0.9,0)--(6,6)--(10,10)--(10,0)--cycle;
 \node at (8,2) {$\mathcal{D}_1$};
 \draw[white,opacity=0.5,fill=gray,rotate around={45:(4.9,4.9)}] (4.9,4.9) ellipse (0.98 and 0.84);

\begin{scope}
\clip (-10,6) rectangle (10,10);
\draw (3,5) .. controls (3.7,3.3) .. (10,10);
\draw (3,-5) .. controls (3.7,-3.3) .. (10,-10);
\draw (-3,5) .. controls (-3.7,3.3) .. (-10,10);
\draw (-3,-5) .. controls (-3.7,-3.3) .. (-10,-10);
\end{scope}
\begin{scope}
\clip (-10,-6) rectangle (10,-10);
\draw (3,5) .. controls (3.7,3.3) .. (10,10);
\draw (3,-5) .. controls (3.7,-3.3) .. (10,-10);
\draw (-3,5) .. controls (-3.7,3.3) .. (-10,10);
\draw (-3,-5) .. controls (-3.7,-3.3) .. (-10,-10);
\end{scope}

\fill (3,5) circle (0.15);
\fill (3,-5) circle (0.15);
\fill (-3,5) circle (0.15);
\fill (-3,-5) circle (0.15);
\end{tikzpicture}
\caption{A possible choice of branch cuts in the case $-\frac{\alpha^2}2<b<\big(2+\sqrt{6}\big)\alpha^2$.
 The branch points of~$\Omega$ are marked with a dot and the branch cuts are represented by dotted lines. The set~$U$ is shaded dark gray and the set $\mathcal{D}_1$ is shaded light gray.}\label{figure3-1}
\end{figure}

In this case $\Omega$ has a branch point in each quadrant of the complex plane.
The contour $\im \Omega(k)=0$, shown in Fig.~\ref{figure3}, consists of the coordinate axes together with two simple curves connecting each of the pairs of zeroes $\big\{\frac{-\alpha + {\rm i}\sqrt{\alpha^2+2 b}}2,\frac{\alpha + {\rm i}\sqrt{\alpha^2+2 b}}2\big\}$ and $\big\{\frac{-\alpha - {\rm i}\sqrt{\alpha^2+2 b}}2,\frac{\alpha - {\rm i}\sqrt{\alpha^2+2 b}}2\big\}$ and intersecting the imaginary axis at
\begin{gather*}
\frac{1}{2} \sqrt{b+\sqrt{-2 \alpha ^4+b^2-4 \alpha ^2 b}}{\rm i}\qquad\text{and}\qquad-\frac{1}{2} \sqrt{b+\sqrt{-2 \alpha ^4+b^2-4 \alpha ^2 b}}{\rm i},
\end{gather*}
 respectively, as well as two parabola like curves intersecting the imaginary axis at
 \begin{gather*}
 \frac{1}{2} \sqrt{b-\sqrt{-2 \alpha ^4+b^2-4 \alpha ^2 b}}{\rm i}\qquad\text{and}\qquad-\frac{1}{2} \sqrt{b-\sqrt{-2 \alpha ^4+b^2-4 \alpha ^2 b}}{\rm i},
 \end{gather*}
 asymptoting towards the lines ${\rm e}^{{\rm i}\frac{\pi}4}$ and ${\rm e}^{{\rm i}\pi- {\rm i}\frac{\pi}4}$, and ${\rm e}^{{\rm i}\pi + {\rm i}\frac{\pi}4}$ and ${\rm e}^{- {\rm i}\frac{\pi}4}$, respectively.

Since there are no branch cuts in $\bar{{D}}_1$, we cannot rule out the corresponding triples using Lemma~\ref{inadmissiblelemma}.
This leads to the following family of potentially admissible triples
\begin{gather*}
\left\{\left(\alpha,\,\omega=\frac{\alpha ^4}{2}-b^2+\alpha ^2 b,\,c=\alpha b{\rm i}\right)\bigg|\,\big(2+\sqrt{6}\big)\alpha^2\leq b,\,\alpha>0\right\}
\end{gather*}
or
\begin{gather*}
\left\{\left(\alpha,\omega,c=\alpha \left(\frac{\alpha ^2}2+\sqrt{\frac{3 \alpha ^4}4- \omega }\right){\rm i}\right)\Bigg|\,\alpha>0,\omega\leq-\frac{1}{2} \big(6 \sqrt{6} +15 \big)\alpha^4\right\}.
\end{gather*}

\vspace{-2mm}

\subsection*{Acknowledgements} The author thanks Jonatan Lenells for helpful discussions.
The author also thanks the anonymous referees for many helpful suggestions.
Support is acknowledged from the European Research Council, Grant Agreement No.~682537.

\vspace{-2mm}

\pdfbookmark[1]{References}{ref}
\LastPageEnding

\end{document}